\documentclass[leqno,12pt]{article}
\usepackage{latexsym}
\usepackage{amsmath}
\usepackage{amssymb}
\usepackage[authoryear]{natbib}
\usepackage{graphicx}
\usepackage{epsfig}

 0

\newcommand{\C}{\mathbb{C}}
\newcommand{\R}{\mathbb{R}}
\newcommand{\D}{\partial\mathbb{D}}
\newcommand{\N}{\mathbb{N}}
\newcommand{\bea}{\begin{eqnarray}}
\newcommand{\eea}{\end{eqnarray}}
\newcommand{\be}{\begin{equation}}
\newcommand{\ee}{\end{equation}}
\newcommand{\beo}{\begin{equation*}}
\newcommand{\eeo}{\end{equation*}}
\newcommand{\beao}{\begin{eqnarray*}}
\newcommand{\eeao}{\end{eqnarray*}}

\newtheorem{theorem}{Theorem}[section]
\newtheorem{lemma}[theorem]{Lemma}

\newtheorem{definition}[theorem]{Definition}

\newtheorem{remark}[theorem]{Remark}
\newtheorem{corollary}[theorem]{Corollary}

\renewcommand{\theequation}{\thesection.\arabic{equation}}
\makeatletter\@addtoreset{equation}{section}\makeatother

\newcommand{\ba}{\begin{array}}
\newcommand{\ea}{\end{array}}
\newcommand{\beqohne}{\begin{eqnarray*}}
\newcommand{\eeqohne}{\end{eqnarray*}}
\newcommand{\beohne}{\begin{equation*}}
\newcommand{\eeohne}{\end{equation*}}

\newcommand{\la}{\lambda}

\renewcommand{\theequation}{\thesection.\arabic{equation}}

\def\3{\ss}

\newcommand{\beq}{\begin{equation}}
\newcommand{\eeq}{\end{equation}}

\begin{document}

\title{Matrix measures on the unit circle, moment spaces, orthogonal polynomials and the Geronimus relations}

\author{
{\small Holger Dette} \\
{\small Ruhr-Universit\"at Bochum} \\
{\small Fakult\"at f\"ur Mathematik} \\
{\small 44780 Bochum, Germany} \\
{\small e-mail: holger.dette@rub.de}\\
\and
{\small Jens Wagener} \\
{\small Ruhr-Universit\"at Bochum} \\
{\small  Fakult\"at f\"ur Mathematik} \\
{\small 44780 Bochum, Germany} \\
{\small e-mail: jens.wagener@rub.de}\\
}

\maketitle

\begin{abstract}
We study the moment space corresponding to matrix measures on the unit circle. Moment points are characterized by non-negative definiteness of
block Toeplitz matrices. This characterization is used to derive an explicit representation of  orthogonal polynomials  with respect to matrix
measures on the unit circle and to present a geometric definition of  canonical moments. It is demonstrated that these geometrically defined
quantities coincide with the Verblunsky coefficients, which appear in the Szeg\"{o} recursions for the matrix orthogonal polynomials. Finally,
we provide an alternative proof of the Geronimus relations which is based on a simple relation between canonical moments of matrix measures on
the interval [-1,1] and the Verblunsky coefficients corresponding to matrix measures on the unit circle.

\end{abstract}

\medskip



\section{Introduction}
\def\theequation{1.\arabic{equation}}
\setcounter{equation}{0}

In recent years considerable interest has been shown in moment problems, orthogonal polynomials, continued fractions and quadrature formulas
corresponding to matrix measures on the real line or on the unit circle. Early work dates back to \cite{krein1949}, while more recent results on
matrix measures on the real line can be found in the papers of \cite{rodman1990}, \cite{duran1995,duran1996} and \cite{dejolapo2000}   among many
others. Additionally, several authors have discussed matrix measures on the unit circle [see \cite{delgenkam1978}, \cite{geronimo1981},
\cite{marrod1989},
\cite{sinass1994b,sinass1996},
\cite{yakmar2001,yakmar2002}, \cite{cafemove2003}].

The
purpose of the present paper is to investigate some geometric properties of the moment space corresponding to matrix measures on the unit
circle. In Section 2 we present a characterization of the moment space in terms of nonnegative definiteness of block Toeplitz matrices. We
also provide a geometric definition of canonical moments of matrix measures on the unit circle, which generalizes the scalar case discussed
by \cite{dettstud1997} in a nontrivial way. In Section 3 an explicit determinantal representation of orthogonal matrix polynomials with respect
to matrix measures on the unit circle is presented,
which generalizes the classical representation in the one-dimensional case [see e.g.\ \cite{geronimus1962}]. These results are used to
identify the canonical moments as Verblunsky coefficients, which appear in the Szeg\"{o} relations for the corresponding orthonormal
and reversed matrix polynomials [see \cite{delgenkam1978},  \cite{sinass1996} or \cite{dampussim2008}]. In particular our results  provide a
geometric definition of Verblunsky coefficients corresponding to matrix measures on the unit circle. Roughly speaking, the Verblunsky
coefficient of order $m$ can be characterized as the distance of the $m$th trigonometric moment to a center of a matrix disc relative to the diameter
of this disc (see Section 3 for more details). Finally, in Section 4 these results are used to present an alternative proof of the Geronimus relations for monic orthogonal polynomials,
which describe the relation between the coefficients in the three-term recursive relation of orthogonal polynomials with respect to a
matrix measure on a compact interval and the coefficients in the Szeg\"{o} recursion of an associated matrix measure on the unit circle.

\section{The moment space of matrix measure on the unit circle}
\def\theequation{2.\arabic{equation}}
\setcounter{equation}{0}

A matrix measure  $\mu$ on the unit circle
is defined as a $p\times p$ matrix of (real valued) Borel measures $\mu = (\mu_{ij})_{ i,j=1,\dots,p}$
on the unit circle  $\D=\{z\in\C|\ |z|=1\}$ such that for each Borel set $A \subset \D$ the matrix
$\mu (A) $ is nonnegative definite, i.e. $\mu (A) \ge 0$. Throughout this paper we use the usual parametrization
 $z=e^{i\theta},\ \theta\in\ [-\pi,\pi)$ and the notation  $\mu(\theta)$ for the sake of simplicity. The $k$th moment of a matrix measure
 $\mu$  on the unit circle is defined by
\be \label{gamma}
\Gamma_k=\Gamma_k (\mu) =\int_{-\pi}^{\pi}e^{ik\theta}d\mu(\theta)  = \alpha_k + i \beta_k
\quad k \in \mathbb{Z}
\ee
where $\alpha_k=\alpha_k(\mu)  =\int_{-\pi}^{\pi}\cos{(k\theta)}d\mu(\theta),$ $\beta_k=
\beta_k (\mu) \int_{-\pi}^{\pi}\sin{(k\theta)}d\mu(\theta)$ $(k=0,1,\dots)$ are the trigonometric moments and
the dependence on the given measure $\mu$ is omitted in the notation, whenever it is clear from the context.
Throughout this paper let $ m \in \mathbb{N}_0 \ \la (\mu) = (\alpha_0,\alpha_1,\beta_1,\dots,\alpha_m,\beta_m) \in (\mathbb{R}^{p \times p})^{2m+1}$
denote the vector of trigonometric moments of order $m$ and define
\begin{equation}\label{M}
{\cal M}_{2m+1}=\{\la (\mu)  ~|~ \mu\text{ is a matrix measure on  }\D\}\subset(\R^{p\times p})^{2m+1}
\end{equation}
as the $(2m+1)$th moment space of matrix measures on the unit circle. The  set ${\cal M}_{2m+1}$ and its interior Int$({\cal M}_{2m+1})$
can be characterized as follows.

\medskip

\begin{theorem}\label{char2}
$\lambda = (\alpha_0,\alpha_1,\beta_1,\dots,\alpha_m,\beta_m)\in {\cal M}_{2m+1}$ if and only if
\beq \label{charac}
\sum_{i=0}^m\sum_{j=0}^m\mbox{trace}(B_iB_j^{\ast}\Gamma_{i-j})\geq0\quad\forall\ B_0,\dots,B_m\in\C^{p\times p},
\eeq
where the matrices $\Gamma_{-m}, \Gamma_{-m+1}, \dots, \Gamma_m$ are defined in (\ref{gamma}).

$\lambda =(\alpha_0,\alpha_1,\beta_1,\dots,\alpha_m,\beta_m)\in Int({\cal M}_{2m+1})$ if and only if
there is strict inequality in (\ref{charac}) except if
$B_0=\dots=B_m=0$.
\end{theorem}

\medskip

{\bf Proof:} We start with a proof of the first part. Assume  that
 $\lambda \in {\cal M}_{2m+1}$ and consider matrices $B_0,\dots,B_m\in\C^{p\times p}$.
 With the notation
\be
B(\theta)=\sum_{k=0}^mB_ke^{ik\theta}\quad(\theta\in[-\pi,\pi))
\ee
it follows that the polynomial $P(\theta)=B(\theta)(B(\theta))^{\ast}$ is obviously nonnegative definite, i.e.
\be\label{pgr}
P(\theta)=B(\theta)(B(\theta))^{\ast}\geq0\quad\forall\theta\in[-\pi,\pi).
\ee
A straightforward calculation shows that the polynomial $P$ can be represented as
\be
\label{dedarst}
P(\theta)=D_0+\sum_{k=1}^mD_k\cos{(k\theta)}+E_k\sin{(k\theta)},
\ee
where the $p\times p$ matrices  $D_0,\dots,D_m,E_1,\dots,E_m$
are defined by $D_0 = A_0$, and for $k=1,\dots,m$
$$D_k=A_k+A_{-k},\quad E_k=i(A_k-A_{-k})$$
and
$$A_k=\sum_{l=0}^{m-k}B_{k+l}B_l^{\ast}\quad\text{and}\quad A_{-k}=A_k^{\ast}.$$
Because it is easy to see that the moment space ${\cal M}_{2m+1}$ is the convex hull of the set
 $$
 \Bigl\{(aa^{\ast},\cos{(\theta)}aa^{\ast},\sin{(\theta)}aa^{\ast},\dots,\cos{(m\theta)}
 aa^{\ast},\sin{(m\theta)}aa^{\ast})~ \Bigl| ~a\in\C^p,\ \theta\in\ [-\pi,\pi) \Bigr\},$$
a similar argument as in Corollary 2.2 of \cite{detstu2002} now shows
that \eqref{pgr} and \eqref{dedarst} imply
\begin{eqnarray*}
0&\leq&\mbox{trace}(D_0\alpha_0)+\sum_{k=1}^m\mbox{trace}(D_k\alpha_k)+\mbox{trace}(E_k\beta_k)\\
&=&\mbox{trace}\left(\int_{-\pi}^{\pi}d(D_0\mu(\theta))+\sum_{k=1}^m\int_{-\pi}^{\pi}\cos{(k\theta)}d(D_k\mu(\theta))+\int_{-\pi}^{\pi}\sin{(k\theta)}d(E_k\mu(\theta))\right)\\
&=&\mbox{trace}\left(\int_{-\pi}^{\pi}\sum_{k=-m}^me^{ik\theta}d(A_k\mu(\theta))\right)\\
&=&\mbox{trace}\left(\int_{-\pi}^{\pi}\sum_{k=0}^me^{ik\theta}d\left(\sum_{l=0}^{m-k}B_{k+l}B_l^{\ast}\mu(\theta)\right)+\int_{-\pi}^{\pi}\sum_{k=1}^me^{-ik\theta}d\left(\sum_{l=0}^{m-k}B_{l}B_{k+l}^{\ast}\mu(\theta)\right)\right)\\
&=&\mbox{trace}\left(\sum_{k=0}^m\sum_{l=0}^m\int_{-\pi}^{\pi}e^{i(k-l)\theta}d(B_kB_l^{\ast}\mu(\theta))\right)\\
&=&\sum_{k=0}^m\sum_{l=0}^m\mbox{trace}(B_kB_l^{\ast}\Gamma_{k-l}),
\end{eqnarray*}
which proves \eqref{charac}. On the other hand assume that the inequality \eqref{charac}
is satisfied  for all matrices $B_0, \dots, B_m \in \mathbb{C}^{p \times p}$ and consider a nonnegative definite matrix polynomial
\be
P(\theta)=D_0+\sum_{k=1}^mD_k\cos{(k\theta)}+E_k\sin{(k\theta)}\geq0\quad\forall\theta\in[-\pi,\pi).
\ee
with hermitian matrices  $D_0,\dots,D_m,E_1,\dots,E_m \in \mathbb{C}^{p \times p}$.
 It now follows from \cite{malyshev1982}
that there exists a matrix polynomial
$$B(\theta)=\sum_{k=0}^mB_ke^{ik\theta}$$
such that $P(\theta)=B(\theta)(B(\theta))^{\ast}$,
and the same
 calculation as in the first part of the proof  yields
$$\mbox{trace}(D_0\alpha_0)+\sum_{k=1}^m\mbox{trace}(D_k\alpha_k)+\mbox{trace}(E_k\beta_k)=\sum_{i=0}^m\sum_{j=0}^m\mbox{trace}(B_iB_j^{\ast}\Gamma_{i-j})\geq0.$$
By similar arguments as in Lemma 2.3 of \cite{detstu2002} it follows that this is sufficient for $\lambda\in{\cal M}_{2m+1}$ .\\
 Finally, the second part of the Theorem is shown similarly observing the fact that
$(\alpha_0,\alpha_1,\beta_1,\dots,\alpha_m,\beta_m)\in Int({\cal M}_{2m+1})$
if and only if
$$\mbox{trace}(D_0\alpha_0)+\sum_{k=1}^m\mbox{trace}(D_k\alpha_k)+\mbox{trace}(E_k\beta_k)>0$$
for any nonnegative definite polynomial $P(\theta)$ of the form \eqref{dedarst}
with  $P(\theta)\neq 0\ \forall\theta\in[-\pi,\pi)$. This characterization can be shown by the same arguments as
presented in \cite{detstu2002} who proved
a corresponding statement
for the moment space of matrix measures on the interval $[0,1]$. \hfill $\Box$.

\bigskip

Throughout this paper let
\be\label{toepmat}
T_m = T_m(\mu) =\left(
\begin{array}{lll}
\Gamma_0 & \cdots & \Gamma_m\\
\vdots & \ddots & \vdots\\
\Gamma_{-m} & \cdots & \Gamma_0
\end{array}
\right)\in\C^{p(m+1)\times p(m+1)}
\ee
denote the Block Toeplitz matrix, where the blocks $\Gamma_i = \Gamma_i(\mu) \ (i=-m,\dots,m)$ are the moments of a matrix
measure $\mu$ on the unit
circle defined by \eqref{gamma} (note that $T_m$ is hermitian). The following
characterization of the moment space ${\cal M}_{2m+1}$ by nonnegative definiteness of Toeplitz matrices is now easily obtained.

\begin{corollary}\label{korollar2}
Assume that $\la=(\alpha_0,\alpha_1,\beta_1,\dots,\alpha_m,\beta_m) \in (\R^{p \times p})^{2m+1}$ and that
 $T_m$ is defined by  (\ref{toepmat}) with $\Gamma_k =  \alpha_k + i \beta_k$ and
  $\Gamma_{-k} =  \alpha_k - i \beta_k$.
\begin{enumerate}
\item[(a)]$ \la \in {\cal M}_{2m+1}$ if and only if  $T_m\geq0$.
\item[(b)]$ \la \in Int({\cal M}_{2m+1})$ if and only if $T_m>0$.
\end{enumerate}
\end{corollary}

\medskip

{\bf Proof:} We only proof part (a); part (b) is shown by similar arguments.
First assume that $ \la \in {\cal M}_{2m+1}$, then we obtain from Theorem \ref{char2}
for all matrices $B_0,\dots,B_m\in\C^{p\times p}$
$$\sum_{i=0}^m\sum_{j=0}^m\mbox{trace}(B_iB_j^{\ast}\Gamma_{j-i})\geq0 .$$
Consequently, if  $a_0,\dots,a_m\in\C^{p}$, $a=(a_0^T,\dots,a_m^T)^T \in \mathbb{C}^{p(m+1)}$ we put
 $B_i=(a_i,0,\dots,0)\in\C^{p\times p}$ ($i=0,\dots,m$) and it follows
\begin{eqnarray*}
a^{\ast}T_ma&=& \mbox{trace} (aa^{\ast}T_m)
=\sum_{i=0}^m\sum_{j=0}^m\mbox{trace}(a_ia_j^{\ast}\Gamma_{j-i})
=\sum_{i=0}^m\sum_{j=0}^m\mbox{trace}(B_iB_j^{\ast}\Gamma_{j-i})
\geq0,
\end{eqnarray*}
which shows that the matrix  $T_m$ is nonnegative definite. To prove the converse assume that
 $T_m\geq0$, i.e.
\be\label{ungl5}
0\leq a^{\ast}T_ma=\sum_{i=0}^m\sum_{j=0}^m\mbox{trace}(a_ia_j^{\ast}\Gamma_{j-i}).
\ee
for all  $a=(a_0^T,\dots,a_m^T)^T \in \C^{p(m+1)}$. If
$B_0,\dots,B_m\in\C^{p\times p}$, and  $a_j^{(i)}$ denotes the $i$th column of the matrix $B_j$ ($j=0,\ldots ,m$, $i=1,\ldots ,p$), then
$$B_jB_k^{\ast}=\sum_{i=1}^pa_j^{(i)}\left(a_k^{(i)}\right)^{\ast}$$
and we obtain from (\ref{ungl5})
$$\sum_{i=0}^m\sum_{j=0}^m\mbox{trace}(B_iB_j^{\ast}\Gamma_{j-i})=\sum_{k=1}^p\sum_{i=0}^m\sum_{j=0}^m\mbox{trace}\left(a_i^{(k)}\left(a_j^{(k)}\right)^{\ast}
\Gamma_{j-i}\right)\geq0.$$
By Theorem \ref{char2} it follows that $ \la \in {\cal M}_{2m+1}$, which completes the proof of
the Corollary. \hfill $\Box$

\bigskip

With the aid of Theorem \ref{char2} and Corollary \ref{korollar2} we are now able
to define geometrically canonical moments for matrix measures on the unit circle. It turns
out that these geometrically defined quantities are exactly the Verblunsky coefficients of matrix measures on the unit circle
as introduced by \cite{dampussim2008} (see Section 3 where we prove this identity). For this purpose let
$W$ denote a $p\times p$ matrix and define
\be
A = A(W) =\left(
\begin{array}{ccccc}
\Gamma_0 & \Gamma_{1}& \cdots & \Gamma_m & W\\
\Gamma_{-1} & \Gamma_0 &  \cdots & \Gamma_{m-1} & \Gamma_m \\
\vdots  & \vdots & \ddots & \vdots &  \vdots \\
\Gamma_{-m} & \Gamma_{-m+1}  &\cdots & \Gamma_0 & \Gamma_{1}\\
W^{\ast} & \Gamma_{-m} & \cdots & \Gamma_{-1} &\Gamma_0
\end{array}
\right)\in\C^{p(m+2)\times p(m+2)}.
\ee
Let $\Gamma^{(m)} = (\Gamma_{-m},  \Gamma_{-m+1}, \ldots , \Gamma_{m-1}, \Gamma_m) \in (\C^{p \times p})^{2m+1}$
denote a vector of  moments of a matrix measure on the unit circle, that is
$(\alpha_0, \alpha_1, \beta_1, \dots, \alpha_m, \beta_m) \in {\cal M}_{2m+1}$, where $\Gamma_k = \alpha_k + i\beta_k$.
Define $ {\cal P}_{\Gamma^{(m)}}$
as the set of all matrix measures  $\mu$ on the unit circle
with moments of order $m$ given by
 $\Gamma^{(m)}$,
that is $\Gamma_j = \int_{-\pi}^{\pi} e^{ik\theta} d \mu (\theta) $
($j=-m,\ldots ,m$).
By Corollary  \ref{korollar2} it follows that the matrix
 $W$ is the  $(m+2)$th moment of a matrix measure  $\mu \in {\cal P}_{\Gamma^{(m)}}$
 if and only if $A (W) \geq 0.$ We assume without loss of generality
that  $(\alpha_0, \alpha_1, \beta_1, \dots, \alpha_m, \beta_m) \in \mbox{Int}({\cal M}_{2m+1})$ which is equivalent to  $T_m>0$ by Corollary
\ref{korollar2} . From  Theorem 1 in \cite{frikir1987} it follows that
$$A(W)\geq0$$ if and only if there exists a $p\times p$ matrix $U$ with  $UU^{\ast}\leq I_p$ such that the matrix
 $W$ can be represented as
 \bea
W&=&\left(\Gamma_1\dots\Gamma_m\right)T_{m-1}^{-1}\left(\Gamma_{-m}\dots\Gamma_{-1}\right)^{\ast}
+ L_m^{1/2} U R_m^{1/2},
\eea
where the matrices $L_m$ and $R_m$ are defined by
\bea \label{lm}
L_m&=&\Gamma_0-\left(\Gamma_1\dots\Gamma_m\right)T_{m-1}^{-1}
\left(\Gamma_1\dots\Gamma_m\right)^{\ast},\\
R_m&=&\Gamma_0-\left(\Gamma_{-m}\dots\Gamma_{-1}\right)T_{m-1}^{-1}
\left(\Gamma_{-m}\dots\Gamma_{-1}\right)^{\ast} \label{rm},
\eea
respectively.
Note that  the matrices  $L_m$ and $R_m$ are Schur complements of the positive definite matrix $T_m$
and as a consequence are also positive definite [see \cite{hornjohn1985}].
This means that that the matrix $W$ is the  $(m+2)$th moment of the matrix measure
$\mu \in {\cal P}_{\Gamma^{(m)}}$, if and only if it is an element of the ``ball''
\be\label{Matrixball}
K_m:=\left\{W \in\mathbb{C}^{p \times p}| L_m^{-1/2}(W-M_m)R_m^{-1/2}=U, UU^{\ast}\leq I_p \right\},
\ee
where the ``center'' of the ball is given by the matrix
\be \label{center}
M_m=\left(\Gamma_1\dots\Gamma_m\right)T_{m-1}^{-1}\left(\Gamma_{-m}\dots\Gamma_{-1}\right)^{\ast}.
\ee
We are now in a position to define the canonical moments of a matrix measure on the unit circle (or Verblunsky coefficients as shown in
Section 3).

\begin{definition}\label{kanonische Momente}
Let  $\mu$ denote a matrix measure on the unit circle with moments
$\Gamma_k=\alpha_k+i\beta_k$ ($k\geq0$), $ \la_{2m+1} (\mu) =(\alpha_0,\alpha_1,\beta_1,\dots,\alpha_m,\beta_m) \in (\mathbb{R}^{p \times p})^{m+1}
\ (m \geq0)$
and define
\be
N(\mu)=\min\left\{m\in\N \mid  \la_{2m+1} (\mu) \in\partial {\cal M}_{2m+1}\right\},
\ee
as the minimum number $m \in \N$ such that $ \la_{2m+1} $ is a boundary point of the
moment space $ {\cal M}_{2m+1}$ (if $ \la_{2m+1} \in \mbox{Int} ({\cal M}_{2m+1})$ for all
$m \in \N$ we put $N(\mu) =\infty)$. For each  $m=0,\dots,N(\mu)-1$ the quantity
\bea
A_{m+1} = A_{m+1}(\mu) &=&L_m^{-1/2}\left(\Gamma_{m+1}-M_m\right)R_m^{-1/2} \label{can}\\
\nonumber
&=&\left[\Gamma_0-\left(\Gamma_1,\dots,\Gamma_m\right)T_{m-1}^{-1}\left(\Gamma_1,\dots,\Gamma_m\right)^{\ast}\right]^{-1/2}\\\notag
\nonumber
&&\times\left(\Gamma_{m+1}-\left(\Gamma_1,\dots,\Gamma_m\right)T_{m-1}^{-1}\left(\Gamma_{-m},\dots,\Gamma_{-1}\right)^{\ast}\right)\\
\nonumber
&&\times\left[\Gamma_0-\left(\Gamma_{-m},\dots,\Gamma_{-1}\right)T_{m-1}^{-1}\left(\Gamma_{-m},\dots,\Gamma_{-1}\right)^{\ast}\right]^{-1/2}\notag
\eea
is called the $(m+1)$th canonical moment of the matrix measure $\mu$.
\end{definition}

\bigskip

Definition \ref{kanonische Momente} is a generalization of the definition of canonical moments of scalar measures on the unit circle
 in \cite{dettstud1997}. In general the explicit representation of the canonical moments in terms
of the moments $\Gamma_0 , \Gamma_1 , \ldots $ is very difficult. For example
if  $m=0$ we have
\be
A_1=\Gamma_0^{-1/2}\Gamma_1\Gamma_0^{-1/2}
\ee
and in the case  $m=1$ we obtain from Definition \ref{kanonische Momente}
\be
A_2=\left(\Gamma_0-\Gamma_1\Gamma_0^{-1}\Gamma_{-1}\right)^{-1/2}\left(\Gamma_2-\Gamma_1\Gamma_0^{-1}\Gamma_1\right)\left(\Gamma_0-\Gamma_{-1}\Gamma_0\Gamma_1\right)^{-1/2}
\ee
In the following section we will demonstrate that the quantities defined by
Definition \ref{kanonische Momente} are the well known Verblunsky coefficients, which are usually
obtained from the recursive relations of the orthonormal polynomials with respect to
matrix measures on the unit circle [see for example \cite{delgenkam1978} where these matrices do not have any special name, \cite{sinass1996} where they are called reflection coefficients or \cite{dampussim2008}]. For this purpose we use an
explicit determinant representation of the matrix orthogonal polynomials, which is of
own interest and given in the following section.

\section{Orthogonal matrix polynomials }
\def\theequation{3.\arabic{equation}}
\setcounter{equation}{0}

A $p \times  p$ matrix polynomial is a $p\times  p$ matrix with
polynomial entries.  It is of degree $n$ if all the polynomial entries  are of
degree less than or equal to $n$ and is  usually
written in the form
\begin{equation} \label{matrpol}   P(z) = \sum_{i=0}^n A_i z^i.
\end{equation}
with coefficients $A_i~ \in \C^{p \times p}$ and $z \in \mathbb{C}$.
Recall that  for matrix polynomials  $P$ and $Q$
the right and left inner product
are defined by
\bea
\langle P,Q\rangle_R &=&\int_{-\pi}^{\pi}P(e^{i\theta})^{\ast}d\mu(\theta)Q(e^{i\theta}),
\label{right} \\
\langle P,Q\rangle_L&=&\int_{-\pi}^{\pi}P(e^{i\theta})d\mu(\theta)Q(e^{i\theta})^{\ast}, \label{left}
\eea
respectively
[see for example \cite{sinass1996}]. The matrix
polynomials  $P$ and $Q$ are called orthogonal with respect to the right inner product
$\langle \cdot, \cdot \rangle_R$ if
\be
\langle P,Q\rangle_R=0
\ee
and orthogonality with respect to the left inner product $\langle \cdot, \cdot \rangle_L$ is defined analogously.
 The matrix polynomials
 $P_0 (z),P_1 (z),P_2 (z),\dots $ are called orthonormal with respect to the right inner product if for each $m \in \mathbb{N}_0$ $P_m(z)$ is of degree
  $m$,  $P_m(z)$ and  $P_{m^{'}}(z)$ are orthogonal with respect to $\langle \cdot, \cdot \rangle_R$ whenever $m\neq m^{'}$
  and
\be
\langle P_m,P_m\rangle_R=I_p,
\ee
where $I_p$ denotes the $p\times p$ identity matrix. Orthonormal polynomials
with respect to the left inner product  $\langle \cdot, \cdot \rangle_L$ are defined analogously.
Orthonormal polynomials with respect to the inner products $\langle \cdot, \cdot \rangle_R$  and $\langle \cdot, \cdot \rangle_L$
are determined uniquely up to multiplication by unitary matrices. In the following discussion
we will derive an explicit representation of these polynomials in terms of the moments
of  matrix measure $\mu$, which generalizes the well known determinant representation
in the scalar case [see for example \cite{geronimus1946}].

For this purpose consider a matrix measure $\mu$ on the unit circle with  moments
$\Gamma_{-m},\dots,\Gamma_m$ and recall the definition of the corresponding block Toeplitz matrix  $T_m$  in \eqref{toepmat}.
We define for $m\in\N$ matrix polynomials by
\bea\label{Psi}
\Psi_m^R(z) &=&\left(T_{ij}^R(z)\right)_{i,j=1,\dots,p}, \\
 \label{Psi2}
\Psi_m^L(z) &=& \left(T_{ij}^L(z)\right)_{i,j=1,\dots,p},
\eea
where the elements  $T_{ij}^R(z)$ and $T_{ij}^L(z)$ in these matrices are given
by the determinants
\be\label{det1}
T_{ij}^R(z)=\left|
\begin{array}{cccc}
\Gamma_0 & \Gamma_1 & \dots & \Gamma_m\\
\Gamma_{-1} & \Gamma_0 & \dots & \Gamma_{m-1}\\
\vdots & \vdots &  & \vdots\\
\Gamma_{-m+1} & \Gamma_{-m+2} & \dots & \Gamma_{1}\\
\Gamma_{-m}^{ij}(z) & \Gamma_{-m+1}^{ij}(z) & \dots & \Gamma_0^{ij}(z)
\end{array}\right|; \quad i,j=1,\dots,p
\ee
and
\be\label{det2}
T_{ij}^L(z)=\left|
\begin{array}{cccc}
\tilde{\Gamma}_0^{ij}(z) & \Gamma_1 & \dots & \Gamma_m\\
\tilde{\Gamma}_{-1}^{ij}(z) & \Gamma_0 & \dots & \Gamma_{m-1}\\
\vdots & \vdots &  & \vdots\\
\tilde{\Gamma}_{-m}^{ij}(z) & \Gamma_{-m+1} & \dots & \Gamma_{0}\\
\end{array}\right|; \quad i,j=1,\dots,p,
\ee
respectively, and the matrices  $\Gamma_{-m+k}^{ij}$ (and $\tilde{\Gamma}_{-m+k}^{ij}$)
are obtained replacing the  $j$th row (and the $i$th column) in the matrix $\Gamma_{-m+k}$ by  $e_i^Tz^k$ (and $e_jz^{m-k}$).
The following result shows that these polynomials are orthogonal with respect to the given
matrix measure $\mu$.

\begin{theorem}\label{orthpol1}
For a given matrix measure $\mu$ on the unit circle let  $\Psi_m^R(z)$ and $\Psi_m^L(z)$ ($m\in\N$) denote the matrix polynomials defined by
(\ref{Psi}) and (\ref{Psi2}), respectively, then we have
\bea \label{ident1}
\langle z^kI_p,\Psi_m^R\rangle_R &=&0 \ ,   ~~(k=0,\dots,m-1)\ ; \quad \langle z^mI_p,\Psi_m^R\rangle_R=|T_m|I_p \\
&& \nonumber
\\
\langle \Psi_m^L,z^kI_p\rangle_L &=& 0 \ ,  ~~(k=0,\dots,m-1) \ ; \quad \langle \Psi_m^L,z^mI_p\rangle_L=|T_m|I_p. \nonumber
\eea
\end{theorem}

\medskip

{\bf Proof:}
We will only give a proof for the polynomials $\Psi_m^R(z)$, the remaining part
of Theorem \ref{orthpol1} is shown similarly. The element $B_{ij}^R$ in the position $(i,j)$
 of the matrix
$$
B^R:=\langle z^kI,\Psi_m^R\rangle_R=\int_{-\pi}^{\pi}e^{-ik\theta}d\mu(\theta)\left(T_{ij}^R(e^{i\theta})\right)_{i,j=1,\dots,p} \quad (k=0,\dots,m),
$$
is given by
\be\label{Bij}
B_{ij}^R=\sum_{l=1}^p\int_{-\pi}^{\pi}e^{-ik\theta}T_{lj}^R(e^{i\theta})d\mu_{il}(\theta).
\ee
An expansion of the determinant
$T_{lj}^R(e^{i\theta})$ with respect to the $(mp+j)$th row yields
\be\label{ent1}
T_{lj}^R(e^{i\theta})=\sum_{n=0}^m(-1)^{(m+n)p+j+l}e^{in\theta}\left|T_m^{(mp+j),(np+l)}\right|,
\ee
where the matrix  $T_m^{(mp+j),(np+l)}$ is obtained from  $T_m$ by deleting the
$(mp+j)$th row and $(np+l)$th column. If  $\gamma_{n,ij} = \int_{-\pi}^{\pi}e^{in\theta}d\mu_{ij}$ denotes the element  of the matrix
$\Gamma_n$ in the position $(i,j)$, where $n\in\{-m,\dots,m\}$, it follows that
\be\label{ent2}
B_{ij}^R=\sum_{n=0}^m\sum_{l=1}^p(-1)^{(m+n)p+j+l}\left|T_m^{(mp+j),(np+l)}\right|\gamma_{n-k,il}.
\ee
Now it is easy to see that the right hand side of \eqref{ent2} is the determinant of the
matrix $T_m$, where the $(mp+j)$th row has been replaced by the vector
$$
(\gamma_{-k,i1},\dots,\gamma_{-k,ip},\gamma_{-k+1,i1},\dots,\gamma_{-k+1,ip},\dots,\gamma_{m-1-k,i1},\dots,\gamma_{m-1-k,ip},\gamma_{m-k,i1},\dots,\gamma_{m-k,ip})
$$
Consequently, if $k \in \{0,\dots,m-1 \}$ the $(mp+j)$th and $(kp+i)$th row in this matrix
coincide and we have  $B_{ij}^R=0$, which proves the first identity in \eqref{ident1}.\\
For a proof of the second identity we note that in the case  $k=m$ and
$i\neq j$ the same argument yields $B_{ij}=0$. If  $k=m$ and $i=j$  it follows  that
$B_{ij}$ is exactly the determinant of the matrix  $T_m$, which completes the proof of the
first assertion of Theorem \ref{orthpol1}. \hfill $\Box$

\bigskip

In the following discussion we derive several consequences of the representations
(\ref{Psi}) and (\ref{Psi2}), which will be useful to identify the canonical moments
as Verblunsky coefficients. In particular we determine the corresponding leading coefficients
and identify the orthonormal polynomials with respect to the measure $\mu$. For this
purpose recall that a matrix polynomial
of the form \eqref{matrpol} is called monic, if the coefficient of the leading term is
the identity matrix, that is $A_n = I_p$.

\begin{corollary}\label{poly}
For a given matrix measure $\mu$ on the unit circle let $\Psi_m^R(z)$ and $\Psi_m^L(z)$  be defined by  (\ref{Psi}) and  (\ref{Psi2})
and consider for $m \leq N(\mu)$ the matrix polynomials
\bea\label{monischr}
\Phi_m^R(z) &= &\Psi_m^R(z) |T_m|^{-1}R_m, \\
\label{monischl}
\Phi_m^L(z)&=&|T_m|^{-1}  L_m \Psi_m^L(z),
\eea
where the matrices $R_m$ and $L_m$ are defined by \eqref{rm} and \eqref{lm}, respectively.
The polynomials $\Phi_m^R(z)$  (and  $\Phi_m^L(z)  $) are monic orthogonal matrix
polynomials with respect to the right (and  left) inner product $\langle \cdot,\cdot\rangle_R$
(and $\langle \cdot,\cdot\rangle_L$).

Similarly, define for $m \leq N(\mu)$
\bea
\label{orthonormalr}
\phi_m^R(z) &=&\Psi_m^R(z)|T_m|^{-1} R_m^{1/2}, \\
 \label{orthonormall}
\phi_m^L(z)&=&|T_m|^{-1}
L_m^{1/2}\Psi_m^L(z),
\eea
then the matrix polynomial $\phi_m^R(z)$ (and $\phi_m^L(z)$) are
orthonormal polynomials with respect to the right (and left) inner
product  $\langle \cdot,\cdot\rangle_R$ (and  $\langle \cdot,\cdot\rangle_L$). The leading coefficients
of  $\phi_m^R(z)$ and $\phi_m^L(z)$ are given by $R^{-1/2}_m$ and $L^{-1/2}_m$,
respectively.
\end{corollary}

\medskip

{\bf Proof:} In the first part we will prove that the leading coefficients of the polynomials
$\Psi_m^R(z)$ and $\Psi_m^L(z)$ defined by  (\ref{Psi}) and
(\ref{Psi2}) are given by
\bea \label{eq1}
L_m^R &=&|T_m| R_m ^{-1}, \\
L_m^L &=& |T_m|L_m^{-1}, \label{eq2}
\eea
respectively. With these representations we
obtain from  Theorem  \ref{orthpol1}
$$
\langle \Psi_m^R,\Psi_m^R\rangle_R =|T_m|(L_m^R)^{\ast}\ ; \quad  \langle \Psi_m^L,\Psi_m^L\rangle_L =|T_m|(L_m^L)^{\ast},
$$
and the assertion of the Corollary follows by a straightforward calculation.

In order to prove \eqref{eq1}  and \eqref{eq2} we restrict ourselves to the first case;
the second case is shown similarly. Observing the definition of the determinants
$T_{ij}^R(z)$ in  (\ref{det1}) we obtain for the entry in the position
 $(i,j)$  of the leading coefficient of the matrix polynomial
 $\Psi_m^R(z)$
 $$
\left(L_m^R\right)_{ij}=(-1)^{2mp+i+j}|T_m^{(mp+j),(mp+i)}|,
$$
where we have used an expansion of the determinant with respect to the
$(mp+j)$ row and the matrix  $T_m^{(mp+j),(mp+i)}$ is obtained from  $T_m$ by deleting the
$(mp+j)$th row and $(mp+i)$th column. This means that $\left(L_m^R\right)_{ij}$ is the entry in the
position $(mp+i,mp+j)$ of the adjoint of the matrix $T_m$ ($i,j=1,\ldots ,p)$,
and consequently $L^R_m / |T_m |$ is the
 $p\times p$ block in the position $(m+1,m+1)$ of the matrix
  $T_m^{-1}$, which is given by
 $$
 \left(\Gamma_0-(\Gamma_{-m}\dots\Gamma_{-1})T_{m-1}^{-1}(\Gamma_{-m}\dots\Gamma_{-1})^{\ast}\right)^{-1}=R_m ^{-1}
$$
[see e.g.  \cite{hornjohn1985}]. This proves the assertion \eqref{eq1} and completes the proof of the Corollary.
\hfill $\Box$.

\bigskip

We are now in a position to identify the canonical moments introduced in Definition
\ref{kanonische Momente}
as Verblunsky coefficients which are defined as coefficients in the Szeg\"{o} relation of the matrix orthonormal polynomials
$\phi^L_n(z)$ and $\phi^R_n(z)$. For this purpose we introduce
for a given matrix polynomial
$P_n$ of degree $n$
the corresponding
reversed polynomial
$$
\tilde{P}_n(z)=z^nP_n\left(\frac{1}{\overline{z}}\right)^{\ast},
$$
where $\overline{z}$ denotes the complex conjugation of $z \in \mathbb{C}$.
  Obviously we have
for any  $p\times p$ matrix $A$
$$
\widetilde{AP}_n(z)=\tilde{P}_n(z)A^{\ast}.
$$
In the following discussion let  $\kappa_m^R = R^{-1/2}_m $ and
 $\kappa_m^L = L^{-1/2}_m$ ($m=1,\dots,N(\mu)-1$) denote the leading coefficients of the orthonormal
 matrix polynomials  $\phi_m^R(z)$ and  $\phi_m^L(z)$ with respect to the right and left
 inner product induced by the matrix measure $\mu$
 and define the matrices
\be \label{rhos}
\rho_m^R=\left(\kappa_{m+1}^R\right)^{-1}\kappa_m^R \quad\text{and}\quad\rho_m^L=\kappa_m^L\left(\kappa_{m+1}^L\right)^{-1}\quad (m=1,\dots,N(\mu)-1).
\ee
Then it follows from \cite{dampussim2008}
  that there
 exist $p \times p $ matrices $H_m$  such that
  the orthonormal matrix polynomial with respect to the  measure
$\mu$    on the unit circle
 satisfy the Szeg\"o recursions
\bea
z\phi_m^L(z)-\rho_m^L\phi_{m+1}^L(z)=H_{m+1}\tilde{\phi}_m^R(z),\label{szegö1}\\
z\phi_m^R(z)-\phi_{m+1}^R(z)\rho_m^R=\tilde{\phi}_m^L(z)H_{m+1}.\label{szegö2}
\eea
 The matrices $H_m$ are uniquely determined
and called Verblunsky or reflection coefficients, because they were introduced for the scalar case
in two seminal papers by \cite{verblunsky1935,verblunsky1936}.
The final result of this section shows that the Verblunsky
coefficients coincide with the canonical moments introduced
in Definition \ref{kanonische Momente}.

\begin{theorem}
\label{A=H}
Let $\mu$  denote a matrix measure on the unit circle and assume that $0\leq m< N(\mu)$. If $A_{m+1}$ is the $(m+1)$th
canonical moment of $\mu$   defined
in  Definition \ref{kanonische Momente} and $H_{m+1}$ is the $(m+1)$th Verblunsky coefficient
 defined by the Szeg\"o recursions \eqref{szegö1} and  \eqref{szegö2}, then
\be
A_{m+1}=H_{m+1}.
\ee
\end{theorem}

\medskip

{\bf Proof:} Integrating the recursion
(\ref{szegö2}) we obtain
$$
\langle I_p,z\phi_m^R-\phi_{m+1}^R\rho_m^R\rangle_R \ = \  \langle I_p,\tilde{\phi}_m^LH_{m+1}\rangle_R
$$
and
$$
\langle I_p,z\Psi_m^R\rangle_R|T_m|^{-1}R_m^{1/2}\ = \ \langle I_p,\tilde{\Psi}_m^L\rangle_R|T_m|^{-1}L_m^{1/2}H_{m+1},
$$
where we have used the orthogonality of the matrix polynomials
 $\Psi_{m+1}^R(z)$  stated in Theorem  \ref{orthpol1} and the representations
 of the orthonormal  polynomials $\phi_m^R$ and   $\phi_m^L$ in Corollary \ref{poly}.
Observing Theorem \ref{orthpol1}
and the identity
\bea\label{tilde1}
\langle I_p,\tilde{\Psi}_m^L\rangle_R &=&\int_{-\pi}^{\pi}d\mu(\theta)e^{im\theta}\left(\Psi_m^L(e^{i\theta})\right)^{\ast}=\langle z^mI_p,\Psi_m^L\rangle_L
=|T_m| I_p
\eea
 yields
\bea\label{darsth}
H_{m+1} &=& L_m^{-{1/2}}\langle I_p,\tilde{\Psi}_m^L\rangle_R^{-1}\langle I_p,z\Psi_m^R\rangle_RR_m^{1/2} \\
&=& L_m^{-{1/2}} |T_m|^{-1}\langle I_p,z\Psi_m^R\rangle_RR_m^{1/2}. \nonumber
\eea
The matrix polynomial $\Psi_m^R (z)$ has the representation
$$
\Psi_m^R(z)=L_m^Rz^m+\sum_{k=0}^{m-1}K_k^Rz^k,
$$
where  $K_0^R, \ldots , K_{m-1}^R$  denote $p\times p$ matrices and
the leading coefficient $L_m^R$ is given by
\eqref{eq1}. Integrating with respect to $d\mu (\theta) $ gives
$$
\langle I_p,z\Psi_m^R\rangle_R=\langle I_p,z^{m+1}+\sum_{k=0}^{m-1}K_k^R\left(L_m^R\right)^{-1}z^{k+1}\rangle_R|T_m|R_m^{-1},
$$
and it follows from (\ref{darsth}) that
\be\label{darsth2}
H_{m+1}=L_m^{-1/2}\langle I_p,z^{m+1}+\sum_{k=0}^{m-1}K_k^R\left(L_m^R\right)^{-1}z^{k+1}\rangle_RR_m^{-1/2}.
\ee
Observing the definition of the canonical moments in (\ref{can}) and the definition of the center (\ref{center}) the assertion of the Theorem
follows if  the identity
\be\label{beh1}
\langle I_p,z^{m+1}+\sum_{k=0}^{m-1}K_k^R\left(L_m^R\right)^{-1}z^{k+1}\rangle_R \ = \ \Gamma_{m+1}-\left(\Gamma_1\dots\Gamma_m\right)T_{m-1}^{-1}\left(\Gamma_{-m}\dots\Gamma_{-1}\right)^{\ast}.
\ee
can be established.
For this purpose we determine the matrices $K_k^R$ ($k=0,\dots,m-1$) explicitly using the representation
of the orthogonal matrix polynomials $\Psi_m^R(z)$ in (\ref{Psi}). From this definition
it follows
that the element in the  position
$(i,j)$ of the matrix   $K_k^R$ is obtained by deleting the  $(mp+j)$th row and the  $(kp+i)$th column
 in the determinant $T_{ij}^R(z)$ defined by (\ref{det1}), that is
$$
\left(K_k^R\right)_{ij}=(-1)^{(m+k)p+i+j}|T_m^{(mp+j),(kp+i)}|.
$$
Here again $T_m^{(mp+j),(kp+i)}$ denotes the matrix obtained $T_m$  by deleting the
$(mp+j)$th row and  $(kp+i)$th column, which coincides with the entry in the position
$(kp+i,mp+j)$ of the
adjoint of the matrix $T_m$. Consequently, it follows that
$$
\left(K_k^R\right)_{ij}=|T_m|(T_m^{-1})_{kp+i,mp+j},
$$
and the ``vector''
$$\frac{1}{|T_m|}\left(\begin{array}{c}
K_0^R\\
\vdots\\
K_{m-1}^R\end{array}\right)
~\in (\C^{p\times p})^m$$
coincides with the right upper block
of size  $mp\times p$ of the matrix  $T_m^{-1}$. By standard result in linear algebra
this block is given by
\beo
-T_{m-1}^{-1}(\Gamma_{-m}\dots\Gamma_{-1})^{\ast}R_m^{-1},
\eeo
which yields
\beao
\langle I_p,\sum_{k=0}^{m-1}K_k^Rz^{k+1}\rangle_R&=&\sum_{k=0}^{m-1}\Gamma_{k+1}K_k^R\\
&=&(\Gamma_1\dots\Gamma_m)\left((K_0^R)^\ast\dots(K_{m-1}^R)^{\ast}\right)^{\ast}\\
&=&-|T_m|(\Gamma_1\dots\Gamma_m)T_{m-1}^{-1}(\Gamma_{-m}\dots\Gamma_{-1})^{\ast}R_m^{-1}.
\eeao
Combining this result with the identity $
\left(L_m^R\right)^{-1}=R_m|T_m|^{-1}
$
finally gives (\ref{beh1}), which completes the proof Theorem
\ref{A=H}. \hfill $\Box$

\section{Geronimus relations for monic polynomials}
\def\theequation{4.\arabic{equation}}
\setcounter{equation}{0}

In this section we present a new proof of the Geronimus relations, which provide
a representation of the canonical moments (or Verblunsky coefficients) of a symmetric matrix measure
on the unit circle in terms of the coefficients in the recurrence relations
of a sequence of orthogonal polynomials with respect to an associated
matrix measure on the interval $[-1,1]$.  There exists several alternative
proofs of these relations in the literature [see \cite{yakmar2001} and \cite{dampussim2008}], but the one
presented here explicitly uses the theory of  canonical moments of matrix measures
as introduced in \cite{detstu2002}. As a by-product we derive several interesting properties of the Verblunsky coefficients.

To be precise let
 $\mu_C$ denote a symmetric (with respect to the point $0$) matrix measure on the unit disc
(i.e. $\mu_C$ is  invariant with respect to the transformation $\theta\mapsto-\theta$).
 We associate to $\mu_C$ a corresponding matrix measure, say  $\mu_I$, on the the interval $[-1,1]$,
  which is defined by the property
\be\label{szegömap}
\int_{-1}^1f(x)d\mu_I(x)=\int_{-\pi}^{\pi}f(\cos{(\theta)})d\mu_C(\theta)
\ee
for all integrable functions
 $f$ defined on the interval  $[-1,1]$. Note that
the relation $Sz:d\mu_C\mapsto d\mu_I$ is called
 Szeg\"o mapping in the literature,
 where the matrix measure  $\mu_I$ is usually defined on the interval
  $[-2,2]$. We will work with the interval  $[-1,1]$ in this section, because
   this interval is also used in the classical  papers of
  \cite{szego1922} and \cite{geronimus1946} and in the monograph on canonical moments by
      \cite{dettstud1997}.

      Note that the inverse of the Szeg\"{o} mapping
      \eqref{szegömap} is characterized by the property
\be\label{inverseszegö}
\int_{-\pi}^{\pi}g(\theta)d\mu_C(\theta)=\int_{-1}^1g(\arccos{(x)})d\mu_I(x),
\ee
where $g$ denotes any integrable function on $ \D$ with
 $g(\theta)=g(-\theta)$ for all $\theta \in [- \pi, \pi)$. For a proof of the Geronimus relations we need several preparations.
 Our first results shows that the canonical moments (or Verblunsky coefficients)
 of  a symmetric matrix measure on the unit circle are real and symmetric matrices.
 The result was also proved by  \cite{dampussim2008}. We provide
 here an alternative proof, because several steps in the proof are used
 later.

\begin{lemma}\label{Asymm}
For any symmetric matrix measure  $\mu_C$ on the unit circle the corresponding
canonical moments $A_m$ are real and symmetric.
\end{lemma}

\medskip

{\bf Proof:} By the symmetry of the matrix measure
 $\mu_C$ we have $
\Gamma_k=\int_{-\pi}^{\pi}e^{ik\theta}d\mu_C(\theta)=\int_{-\pi}^{\pi}e^{-ik\theta}d\mu_C(\theta)=\Gamma_{-k}
$
which yields $\Gamma_k=\int_{-\pi}^{\pi}\cos{(k\theta)}d\mu_C(\theta)$. Consequently, the
block Toeplitz matrix associated with
$\mu_C$ is given by
\be\label{tsymmetrisch}
T_m=\left(\begin{array}{lll}
\Gamma_0 & \dots & \Gamma_m\\
\vdots & \ddots & \vdots\\
\Gamma_m & \dots & \Gamma_0
\end{array}\right)
\ee
and is symmetric.  Because all entries of the matrix $T_m$ are real,
the canonical moments
 $A_m$ are also real  and it remains to establish the symmetry.\\
 For this purpose we denote by
  $[A]_{(k,l)}$ the  $p\times p$ block in the position
  $(k,l)$ of the  $mp\times mp-$ block matrix $A$. We will show
  at the end of this proof that
  \be \label{inversesymm}
\left[T_{m-1}^{-1}\right]_{(k,l)}=\left[T_{m-1}^{-1}\right]_{(m+1-k,m+1-l)}.
\ee
From this identity and the property $\Gamma_k = \Gamma_k^*$ we obtain
\begin{eqnarray*}
(\Gamma_1,\dots,\Gamma_m)T_{m-1}^{-1}(\Gamma_m,\dots,\Gamma_1)^{\ast}&=&\sum_{k,l=1}^{m}\Gamma_k\left[T_{m-1}^{-1}\right]_{(k,l)}\Gamma_{m+1-l}
=
\sum_{k,l=1}^{m}\Gamma_{m-k+1}\left[T_{m-1}^{-1}\right]_{(m-k+1,m-l+1)}\Gamma_{l}\\
&=&\sum_{k,l=1}^{m}\Gamma_{m-k+1}\left[T_{m-1}^{-1}\right]_{(k,l)}\Gamma_{l}
= (\Gamma_m,\dots,\Gamma_1)T_{m-1}^{-1}(\Gamma_1,\dots,\Gamma_m)^{\ast},
\end{eqnarray*}
and by similar arguments
\be\label{L=R}
\left(\Gamma_1,\dots,\Gamma_m\right)T_{m-1}^{-1}\left(\Gamma_1,\dots,\Gamma_m\right)^{\ast}=\left(\Gamma_{m},\dots,\Gamma_{1}\right)
T_{m-1}^{-1}\left(\Gamma_{m},\dots,\Gamma_{1}\right)^{\ast}.
\ee
Observing the definition of the canonical moments $A_{m+1}$
it now follows that
\beao
A_{m+1}^{\ast}&=&\left[\Gamma_0-\left(\Gamma_{m},\dots,\Gamma_{1}\right)T_{m-1}^{-1}\left(\Gamma_{m},\dots,\Gamma_{1}
\right)^{\ast}\right]^{-1/2}\left(\Gamma_{m+1}-\left(\Gamma_1,\dots,\Gamma_m\right)T_{m-1}^{-1}\left(\Gamma_{m},\dots,\Gamma_{1}
\right)^{\ast}\right)^{\ast}\\
&&\times\left[\Gamma_0-\left(\Gamma_1,\dots,\Gamma_m\right)T_{m-1}^{-1}\left(\Gamma_1,\dots,\Gamma_m\right)^{\ast}\right]^{-1/2}\\
&=&A_{m+1}
\eeao
which proves the remaining assertion of Lemma \ref{Asymm}.

\bigskip

{\it Proof of the identity \eqref{inversesymm}.} The element in the position
 $(i,j)$ of the matrix $\left[T_{m-1}^{-1}\right]_{(k,l)}$
 and  $\left[T_{m-1}^{-1}\right]_{(m+1-k,m+1-l)}$ are given by
$$|T_{m-1}|^{-1}(-1)^{(l+k)p+i+j}\left|T_{m-1}^{((l-1)p+j),((k-1)p+i)}\right|$$
and
 $$|T_{m-1}|^{-1}(-1)^{(2m-l-k)p+i+j}\left|T_{m-1}^{((m-l)p+j),((m-k)p+i)}\right|,$$
 respectively, where $T_{m-1}^{((m-l)p+j),((m-k)p+i)}$ denotes the matrix obtained from
 $T_{m-1}$ by deleting the $(m-l)p+j$ row and $(m-k)p+i$ column
 (note that both expressions have the same sign).
In the following discussion
we denote by  $A^{(\cdot),(i)}$ and   $A^{(j),(\cdot)}$ the matrix obtained from $A$ by deleting the $i$th column or
the $j$th row, respectively. Then interchanging first  columns and then rows
  yields
\beao
\left|T_{m-1}^{((l-1)p+j),((k-1)p+i)}\right|&=&\left|\begin{array}{llrllrr}
\Gamma_0 & \dots & \Gamma_{k-2} & \Gamma_{k-1}^{(\cdot),(i)} & \Gamma_{k} & \dots & \Gamma_{m-1}\\
\vdots & & \vdots & \vdots & \vdots & &\vdots\\
\Gamma_{l-2} & \dots & \Gamma_{|l-k|}& \Gamma_{|l-k-1|}^{(\cdot),(i)} & \Gamma_{|l-k-2|} & \dots & \Gamma_{m-l+1}\\
\Gamma_{l-1}^{(j),(\cdot)} & \dots & \Gamma_{|l-k+1|}^{(j),(\cdot)} & \Gamma_{|l-k|}^{(j),(i)}& \Gamma_{|l-k-1|}^{(j),(\cdot)} & \dots & \Gamma_{m-l}^{(j),(\cdot)}\\
\Gamma_{l} & \dots & \Gamma_{|l-k+2|}& \Gamma_{|l-k+1|}^{(\cdot),(i)} & \Gamma_{|l-k|} & \dots & \Gamma_{m-l-1}\\
\vdots & & \vdots & \vdots & \vdots & &\vdots\\
\Gamma_{m-1} & \dots & \Gamma_{m-k+1} & \Gamma_{m-k}^{(\cdot),(i)} & \Gamma_{m-k-1} & \dots & \Gamma_{0}
\end{array}\right|\\
&=& (-1)^{\gamma}
\left|\begin{array}{llrllrr}
\Gamma_{m-1} & \dots & \Gamma_{k} & \Gamma_{k-1}^{(\cdot),(i)} & \Gamma_{k-2} & \dots & \Gamma_{0}\\
\vdots & & \vdots & \vdots & \vdots & &\vdots\\
\Gamma_{m-l+1} & \dots & \Gamma_{|l-k-2|}& \Gamma_{|l-k-1|}^{(\cdot),(i)} & \Gamma_{|l-k|} & \dots & \Gamma_{l-2}\\
\Gamma_{m-l}^{(j),(\cdot)} & \dots & \Gamma_{|l-k-1|}^{(j),(\cdot)} & \Gamma_{|l-k|}^{(j),(i)}& \Gamma_{|l-k+1|}^{(j),(\cdot)} & \dots & \Gamma_{l-1}^{(j),(\cdot)}\\
\Gamma_{m-l-1} & \dots & \Gamma_{|l-k|}& \Gamma_{|l-k+1|}^{(\cdot),(i)} & \Gamma_{|l-k+2|} & \dots & \Gamma_{l}\\
\vdots & & \vdots & \vdots & \vdots & &\vdots\\
\Gamma_{0} & \dots & \Gamma_{m-k-1} & \Gamma_{m-k}^{(\cdot),(i)} & \Gamma_{m-k+1} & \dots & \Gamma_{m-1}
\end{array}\right|\\
&=&(-1)^{2\gamma}
\left|\begin{array}{llrllrr}
\Gamma_{0} & \dots & \Gamma_{m-k-1} & \Gamma_{m-k}^{(\cdot),(i)} & \Gamma_{m-k+1} & \dots & \Gamma_{m-1}\\
\vdots & & \vdots & \vdots & \vdots & &\vdots\\
\Gamma_{m-l-1} & \dots & \Gamma_{|l-k|}& \Gamma_{|l-k+1|}^{(\cdot),(i)} & \Gamma_{|l-k+2|} & \dots & \Gamma_{l}\\
\Gamma_{m-l}^{(j),(\cdot)} & \dots & \Gamma_{|l-k-1|}^{(j),(\cdot)} & \Gamma_{|l-k|}^{(j),(i)}& \Gamma_{|l-k+1|}^{(j),(\cdot)} & \dots & \Gamma_{l-1}^{(j),(\cdot)}\\
\Gamma_{m-l+1} & \dots & \Gamma_{|l-k-2|}& \Gamma_{|l-k-1|}^{(\cdot),(i)} & \Gamma_{|l-k|} & \dots & \Gamma_{l-2}\\
\vdots & & \vdots & \vdots & \vdots & &\vdots\\
\Gamma_{m-1} & \dots & \Gamma_{k} & \Gamma_{k-1}^{(\cdot),(i)} & \Gamma_{k-2} & \dots & \Gamma_{0}
\end{array}\right|\\
&& \\
 \\
&=&\left|T_{m-1}^{((m-l)p+j),((m-k)p+i)}\right|,
\eeao
for some $\gamma \in \mathbb{N}$, because the number of changed columns coincides with the number of changed rows.
This implies (\ref{inversesymm}) an completes the proof of Lemma \ref{Asymm}. \hfill $\Box$

\bigskip

For the next step we need to define canonical moments of matrix measures on the
interval $[-1,1]$. Because the main arguments here are very similar to the proceeding in
\cite{detstu2002}, who considered  matrix measures on the interval $[0,1]$,
we only state the main differences without proofs. To be precise, define for a matrix measure
$\mu_I$ on the interval
$[-1,1]$ the moments $S_k = S_k ( \mu_I)
= \int_{-1}^1 x^k d\mu_I (x)$ ($k=0,1,\ldots $)
and a vector $c_n( \mu_I) = ( S_0 ( \mu_I), \ldots , S_n ( \mu_I)) \in (\R^{p\times p})^{n+1}$.
We consider the moment space
\be \label{momspaceint}
 {\cal M}_{n+1}^{(I)} =
 \{c_n (\mu_I)  ~|~ \mu_I\text{ is a matrix measure on  } [-1,1] \}\subset(\R^{p\times p})^{n+1} \ee
 corresponding to the first $n$ moments of matrix measures on the interval $[-1,1]$. For  a matrix measure $\mu_I$ on the interval $[-1,1]$
we define the
block Hankel matrices $\overline{H}_j$ and $\underline{H}_j$
\beao
\underline{H}_{2m}&=&\left(
\begin{array}{ccc}
S_0 & \dots & S_m\\
\vdots & \ddots  & \vdots\\
S_m & \dots & S_{2m}
\end{array}\right),
\\
\overline{H}_{2m}&=&\left(
\begin{array}{ccc}
S_0-S_2 & \dots & S_{m-1}-S_{m+1}\\
\vdots &  \ddots & \vdots\\
S_{m-1}-S_{m+1} & \dots & S_{2m-2}-S_{2m}
\end{array}\right),
\\
\underline{H}_{2m+1}&=&\left(
\begin{array}{ccc}
S_0+S_1 & \dots & S_m+S_{m+1}\\
\vdots & \ddots  & \vdots\\
S_m+S_{m+1} & \dots & S_{2m}+S_{2m+1}
\end{array}\right),
\\
\overline{H}_{2m+1} &=& \left(
\begin{array}{ccc}
S_0-S_1 & \dots & S_m-S_{m+1}\\
\vdots & \ddots  & \vdots\\
S_m-S_{m+1} & \dots & S_{2m}-S_{2m+1}
\end{array}\right).
\eeao
We introduce the notation
\beao
\underline{h}_{2m}&=&(S_m,\dots, S_{2m-1})^T,\quad \overline{h}_{2m}=(S_{m-1}-S_{m+1},\dots,S_{2m-3}-S_{2m-1})^T ,\\
\underline{h}_{2m+1}&=&(S_m+S_{m+1},\dots,S_{2m-1}+S_{2m})^T,\quad\overline{h}_{2m+1}=(S_m-S_{m+1},\dots,S_{2m-1}-S_{2m})^T,
\eeao
and define  $S_1^{+}=S_0$, $S_2^{+}=S_0$,
\bea \nonumber
S_{2m}^{+} &=&S_{2m-2}-\overline{h}_{2m}^T\overline{H}_{2m-2}^{-1}\overline{h}_{2m}\quad(m\geq2), \\
&& \label{splus} \\
S_{2m+1}^{+} &=& S_{2m}-\overline{h}_{2m+1}^T\overline{H}_{2m-1}^{-1}\overline{h}_{2m+1}\quad(m\geq1),
\nonumber
\eea
and  $S_1^{-}=-S_0$,
\bea
S_{2m}^{-}&=&\underline{h}_{2m}^T\underline{H}_{2m-2}^{-1}\underline{h}_{2m}\quad(m\geq1),\nonumber \\
&& \label{sminus} \\
S_{2m+1}^{-} &=&\underline{h}_{2m+1}^T\underline{H}_{2m-1}^{-1}\underline{h}_{2m+1}-S_{2m}\quad(m\geq1).
\nonumber
\eea
Note that the quantities $S_n^+$ and $S_n^-$ are determined by  $S_0,\ldots , S_{n-1}$.
It can be shown by the same argument as in \cite{detstu2002} that for $(S_0,\dots,S_{n-1})\in \mbox {Int} ({\cal M}_n)$ and any matrix measure
$\mu_I$
on the interval $[-1,1]$ with moments satisfying $S_j(\mu_I)=S_j$ $(j=0,\dots,n-1)$,  the  moment of order $n$
$S_n ( \mu_I) = \int_{-1}^1 x^n d \mu_I (x)$ satisfies
\be\label{schranken}
S_{n}^{-}\leq S_{n} ( \mu_I) \leq S_{n}^{+},
\ee
With these preparations we can define the canonical moments
of a matrix measure on the interval $[-1,1]$ with moments  $S_0,\dots,S_{n-1}$.

\begin{definition}\label{intervall}
Let
$\mu_I$ denote a matrix measure on the interval $[-1,1]$
with moments $S_k = S_k ( \mu_I)
= \int_{-1}^1 x^k d\mu_I (x)$ ($k=0,1,\ldots $) and define
\be
N(\mu_I)=\min\left\{k\in\N \mid (S_0,\dots,S_k)\in\partial {\cal M}_{k+1}^{(I)}\right\}.
\ee
For any $n=0,\dots,N(\mu_I)-1$ the
(symmetric)    canonical moments of the matrix measure $\mu_I$ are defined by
\be
U_{n+1}=\left(S_{n+1}^{+}-S_{n+1}^{-}\right)^{-1/2}\left(S_{n+1}-S_{n+1}^{-}\right)\left(S_{n+1}^{+}-S_{n+1}^{-}\right)^{-1/2},
\ee
where the quantities
 $S_{n+1}^{+}$ and $S_{n+1}^{-}$ are given by (\ref{splus}) and (\ref{sminus}), respectively.
\end{definition}

\bigskip

Note that \cite{detstu2002} use a non symmetric definition of canonical moments of matrix measures on the interval $[0,1]$, that
is
\be\label{nonsym}
\bar
U_{n+1}=\left(S_{n+1}^{+}-S_{n+1}^{-}\right)^{-}\left(S_{n+1}-S_{n+1}^{-}
\right).
\ee
This non symmetric definition turns out to be more useful when working with monic orthogonal polynomials but in the present context the symmetric version has advantages.
We are now in a position to prove the main result of this section, which relates the
canonical moments
of a symmetric matrix measure on the unit circle and the canonical
moments of the associated matrix measure on the interval $[-1,1]$ by
the Szeg\"{o} mapping.
For this purpose recall  the definition of the matrix ball $K_m$
 in (\ref{Matrixball}) and the defintion for the matrices  $L_m$, $R_m$ and  $M_m$
 \eqref{lm}, \eqref{rm} and \eqref{center}, respectively. If the given
 measure  $\mu_C$  on the  unit circle is symmetric, then it follows from
(\ref{L=R})
\be \label{4.12}
L_m=R_m.
\ee
The following result is the main step for the proof of the Geronimus relations.

\begin{theorem} \label{geron1}
Let   $\mu_C$ denote  a symmetric matrix measure on the unit circle and denote by $\mu_I = Sz(\mu_C)$ the associated
matrix measure on the interval $[-1,1]$ defined by the Szeg\"o mapping  (\ref{szegömap}).
The canonical moments  $A_n$ and  $U_n$ of the matrix measures $\mu_C$ and
 $\mu_I$ satisfy
$$
A_n=2U_n-I_p \ ;\quad n=1,\dots,N(\mu_C).
$$
Similarly,
the non symmetric canonical moments $\overline{U}_n$ defined in (\ref{nonsym})
satisfy
\be \label{d1}
2\overline{U}_n-I_p=\overline{A}_n \ ;\quad n=1,\dots,N(\mu_C),
\ee
where the quantities $\overline{A}_n$ are  given by
\be \label{d2}
\overline{A}_n=L_{n-1}^{-1/2}A_nL_{n-1}^{1/2}.
\ee
\end{theorem}

\medskip

{\bf Proof:}
We only prove the first part of the Theorem. The second part is shown
by similar
arguments.
Assume that  $m<N(\mu_C)$ and let $\Gamma_0,\Gamma_1, \ldots,$ denote  moments
of the matrix measure on the unit circle $\mu_C$. For $j=0,1, \ldots $ we define $T_j (x) = \cos (j \arccos x) $ as the $j$th
 (scalar) Chebychev polynomial of the first kind, then it follows from (\ref{inverseszegö}) and from \cite{rivlin1990} that
\bea\label{chebychev}
\Gamma_j&=&\int_{-\pi}^{\pi}\cos{(j\theta)}d\mu_C(\theta)=\int_{-1}^1T_j(x)d\mu_I(x)\notag\\
&=&\sum_{k=0}^{\lfloor j/2\rfloor}(-1)^k\frac{j\Gamma(j-k)}{\Gamma(k+1)\Gamma(j-2k+1)}2^{j-2k-1}S_{j-2k},
\eea
where $S_l=\int_{-1}^1x^ld\mu_I(x)$ $(l=0,1,\dots)$ denote the moments of the associated matrix measure
$\mu_I=Sz(\mu_C)$ on the interval. Recall the definition of  $S_{m+1}^{+}$ and $S_{m+1}^{-}$ in
     (\ref{splus}) and (\ref{sminus}), then there  exist matrix measures
      $\mu_I^{+}$ and $\mu_I^{-}$ on the interval  $[-1,1]$
      such that $S_j = S_j (\mu^\pm_I)$ ($j=0,\dots ,m$) and
$$
S_{m+1}^{+}=\int_{-1}^1x^{m+1}d\mu_I^{+}(x)\quad\text{and}\quad S_{m+1}^{-}=\int_{-1}^1x^{m+1}d\mu_I^{-}(x).
$$
We define
\bea
\Gamma_{m+1}^{+}& = &2^mS_{m+1}^{+}+\sum_{k=1}^{\lfloor (m+1)/2\rfloor}(-1)^k\frac{(m+1)\Gamma(m+1-k)}{\Gamma(k+1)\Gamma(m-2k+2)}2^{m-2k}S_{m+1-2k} \label{Gammaplus}\\
\Gamma_{m+1}^{-} &= &2^mS_{m+1}^{-}+\sum_{k=1}^{\lfloor (m+1)/2\rfloor}(-1)^k\frac{(m+1)\Gamma(m+1-k)}{\Gamma(k+1)\Gamma(m-2k+2)}2^{m-2k}S_{m+1-2k}\label{Gammaminus}.
\eea
With the inverse Szeg\"o mapping we obtain the symmetric measures
$ \mu_C^{+} =  (Sz)^{-1}(\mu_I^{+})$ and $\mu_C^{-} =
(Sz)^{-1}(\mu_I^{-})$ on the unit circle and the representation (\ref{chebychev}) yields that
 the  measures $\mu_C^{-}$ and $\mu_C^{+}$ satisfy
$$
\int_{-\pi}^{\pi}\cos{((m+1)\theta)}d\mu_C^{+}(\theta)=\Gamma_{m+1}^{+}\quad\text{and}\quad\int_{-\pi}^{\pi}\cos{((m+1)\theta)}d\mu_C^{-}(\theta)=\Gamma_{m+1}^{-}.
$$
Consequently, recalling the definition of the set $K_m$ in (\ref{Matrixball}) we have
$\Gamma_{m+1}^{+}, \Gamma_{m+1}^{-} \in K_m$ and from the extremal property of the
moments $S_{m+1}^{+}$ and  $S_{m+1}^{-}$ we obtain that $\Gamma_{m+1}^{+}, \Gamma_{m+1}^{-} \in \partial K_m$.
 By the definition of the set $K_m$ in \eqref{Matrixball} it therefore follows that
the canonical moments $A_{m+1}^{+}$ and  $A_{m+1}^{-}$ corresponding to matrix measures
$\mu_C^{+}$ and  $\mu_C^{-}$, respectively, are unitary.
Moreover, Lemma \ref{Asymm}, implies that the matrices  $A_{m+1}^{+}$ and $A_{m+1}^{-}$ are symmetric
 with real entries, which yields
$$
\left(A_{m+1}^{+}\right)^2=I_p\quad\text{and}\quad\left(A_{m+1}^{-}\right)^2=I_p.
$$
Consequently all eigenvalues of the matrices  $A_{m+1}^{+}$ and $A_{m+1}^{-}$
are given by  $-1$ and $1$. \

We now define the matrices
\be \label{erg1}
\tilde{\Gamma}_{m+1}^{+}=M_m+L_m\quad\text{and}\quad\tilde{\Gamma}_{m+1}^{-}=M_m-L_m,
\ee
which are obviously elements of the set $K_m$ because by (\ref{4.12}) we have $L_m=R_m$. Consequently, there exist matrix measures
 $ \tilde{\mu}_C^{+}$ and  $\tilde{\mu}_C^{-}$ such that
 $\Gamma_j ( \tilde{\mu}_C^{\pm}) = \Gamma_j$ ($j=0,\ldots ,m$) and
 \beao
 \Gamma_{m+1} ( \tilde{\mu}_C^{+}) &=& \tilde{\Gamma}_{m+1}^{+}\\
\Gamma_{m+1} ( \tilde{\mu}_C^{-}) &=& \tilde{\Gamma}_{m+1}^{-}
\eeao
Without loss of generality we assume that $\tilde{\mu}_C^{+}$ and $\tilde{\mu}_C^{-}$ are
symmetric with respect to the point $0$ [otherwise use $\frac{1}{2}(\tilde{\mu}_C^{+}(\theta)+\tilde{\mu}_C^{+}(-\theta))$] and we define
$\tilde{\mu}_I^+ = Sz (\tilde{\mu}_C^{+}) $ and $\tilde{\mu}_I^- = Sz (\tilde{\mu}_C^{-}) $
as the associated measures on the interval $[-1,1]$ with $(m+1)$th moments $\tilde{S}_{m+1}^{+}$ and
 $\tilde{S}_{m+1}^{-}$, respectively. These matrices satisfy the identities
\beao
\tilde{\Gamma}_{m+1}^{+}&=&2^m\tilde{S}_{m+1}^{+}+\sum_{k=1}^{\lfloor (m+1)/2\rfloor}(-1)^k\frac{(m+1)\Gamma(m+1-k)}{\Gamma(k+1)\Gamma(m-2k+2)}2^{m-2k}S_{m+1-2k} \\
\tilde{\Gamma}_{m+1}^{-}&=&2^m\tilde{S}_{m+1}^{-}+\sum_{k=1}^{\lfloor (m+1)/2\rfloor}(-1)^k\frac{(m+1)\Gamma(m+1-k)}{\Gamma(k+1)\Gamma(m-2k+2)}2^{m-2k}S_{m+1-2k}
\eeao
From the inequalities (\ref{schranken}) it follows that $S_{m+1}^{+}\geq \tilde{S}_{m+1}^{+}$ and $\tilde{S}_{m+1}^{-}\geq S_{m+1}^{-}$
 (note that $\tilde{S}_{m+1}^{+}$ and $\tilde{S}_{m+1}^{-}$ are moments of a matrix measure on the interval
$[-1,1]$ with moments $S_0,\dots,S_m$). On the other hand we have
\beao
2^m\left(\tilde{S}_{m+1}^{+}-S_{m+1}^{+}\right)&=&\tilde{\Gamma}_{m+1}^{+}-\Gamma_{m+1}^{+}\\
&=&M_m+L_m-(M_m+L_m^{1/2}A_{m+1}^{+}L_m^{1/2})\\
&=&L_m^{1/2}\left(I_p-A_{m+1}^{+}\right)L_m^{1/2}\\
&\geq&0,
\eeao
because the eigenvalues of the
matrix  $I_p-A_{m+1}$ are given by
 $0$ and $2$. So we obtain
$$
\tilde{S}_{m+1}^{+}=S_{m+1}^{+},
$$
while a similar argument shows
$$
\tilde{S}_{m+1}^{-}=S_{m+1}^{-}.
$$
Consequently, it follows that
\beao
A_{m+1}^{+}&=& I_p \ ; \quad \quad A_{m+1}^{-}=-I_p \ ; \\
\tilde{\Gamma}_{m+1}^{+}
&=&\Gamma_{m+1}^{+} \ ; \quad\quad\tilde{\Gamma}_{m+1}^{-}=\Gamma_{m+1}^{-} \ ;
\eeao
and we obtain from the definitions of  $\tilde{\Gamma}_{m+1}^{+}$,
$\tilde{\Gamma}_{m+1}^{-}$ in (\ref{erg1})
$$
M_m=\frac{1}{2}(\Gamma_{m+1}^++\Gamma_{m+1}^-),\quad L_m=\frac{1}{2}(\Gamma_{m+1}^+-\Gamma_{m+1}^-).
$$
The definition of the ($m+1$)th canonical moment $A_{m+1}$ of the matrix measure $\mu$ and (\ref{Gammaplus})-(\ref{Gammaminus}) now imply
\beao
A_{m+1}&=&L_m^{-1/2}(\Gamma_{m+1}-M_m)L_m^{-1/2}\\
&=&\Bigl(\frac{1}{2}\bigl(\Gamma_{m+1}^{+}-\Gamma_{m+1}^{-}\bigr)\Bigr)^{-1/2}\Bigl(\Gamma_{m+1}
-\frac{1}{2}\bigl(\Gamma_{m+1}^{+}+\Gamma_{m+1}^{-}\bigr)\Bigr)\Bigl(\frac{1}{2}
\bigl(\Gamma_{m+1}^{+}-\Gamma_{m+1}^{-}\bigr)\Bigr)^{-1/2}\\
&=&\left(S_{m+1}^{+}-S_{m+1}^{-}\right)^{-1/2}\left(2S_{m+1}-(S_{m+1}^{+}+S_{m+1}^{-})\right)\left(S_{m+1}^{+}-S_{m+1}^{-}\right)^{-1/2}\\
&=&2\left(S_{m+1}^{+}-S_{m+1}^{-}\right)^{-1/2}\left(S_{m+1}-S_{m+1}^{-}\right)\left(S_{m+1}^{+}-S_{m+1}^{-}\right)^{-1/2}-I_p\\
&=&2U_{m+1}-I_p,
\eeao
where the last equality
is a consequence of  the definition of canonical moments of matrix measures on the interval $[-1,1]$. This proves the assertion of the theorem.
\hfill $\Box$

\bigskip

Our final result gives the Geronimus relations for monic orthogonal matrix polynomials, which generalize
the results obtained by \cite{geronimus1946} and \cite{faybusovich} for the scalar case.
To be precise note that  Corollary 3.2
together with (\ref{4.12}) yield for the monic orthogonal polynomials $\Phi_m^R$ and $\Phi_m^L$ defined in (\ref{monischr}) and (\ref{monischl}), respectively
\beao
\rho_m^L\phi_{m+1}^L&=&L_m^{-1/2}\Phi_{m+1}^L,\quad \phi_{m+1}^R\rho_m^R=\Phi_{m+1}^RL_m^{-1/2} \\
\tilde{\phi}_m^R&=&L_m^{-1/2}\tilde{\Phi}_m^R,\quad \tilde{\phi}_m^L=\tilde{\Phi}_m^LL_m^{-1/2}.
\eeao
Using these equations we obtain from (\ref{szegö1}), (\ref{szegö2}) and the second part of Theorem \ref{geron1} the Szeg\"o recursion for the monic orthogonal matrix polynomials with
respect to a matrix measure on the unit circle, that is
\beao
z\Phi_m^L(z)-\Phi_{m+1}^L(z)&=&\overline{A}_{m+1}^{\ast}\tilde{\Phi}_m^R(z),\\
z\Phi_m^R(z)-\Phi_{m+1}^R(z)&=&\tilde{\Phi}_m^R(z)\overline{A}_{m+1}
\eeao
Consequently, the matrices $\overline{A}_{m+1}$ defined by
\eqref{d2}
are the Verblunsky coefficients corresponding to the  monic orthogonal polynomials and we obtain the following result.

\bigskip

\begin{theorem}
Let   $\mu_C$ denote  a symmetric matrix measure on the unit circle and denote by $\mu_I = Sz(\mu_C)$ the associated
matrix measure on the interval $[-1,1]$ defined by the Szeg\"o mapping  (\ref{szegömap}).
If  $P_0$, $P_1$,\dots be the monic polynomials orthogonal with respect to the matrix
measure $\mu_I$ satisfying the three term recurrence recursion
\be  \label{dreischritt}
(1+t)P_{m+1}(t)=P_{m+2}(t)+P_{m+1}(t)C_{m+1}+P_m(t)B_m,
\ee
($P_0(t)=I_p$, $P_{-1}(t)=0_p$), then
the matrices $B_m$ and $C_{m+1}$ satisfy
\beao
B_m &=&\frac{1}{4}(I_p-\overline{A}_{2m})(I_p-\overline{A}_{2m+1}^2)(I_p+\overline{A}_{2m+2}), \\
C_{m+1}&=& \frac{1}{2}(I_p-\overline{A}_{2m+1})
(I_p+\overline{A}_{2m+2})+\frac{1}{2}(I_p-\overline{A}_{2m+2})(I_p+\overline{A}_{2m+3}),
\eeao
where the quantities $\overline{A}_n$ are defined in \eqref{d2}.
\end{theorem}

\medskip

{\bf Proof:} It follows analogously to \cite{detstu2002}  that the
 matrices $B_m$ and $C_{m+1}$ are given by
\beao
B_m &=&(S_{2m}-S_{2m}^{-})^{-1}(S_{2m+2}-S_{2m+2}^-), \\
C_{m+1}&=&(S_{2m+2}-S_{2m+2}^-)^{-1}(S_{2m+3}-S_{2m+3}^-)+(S_{2m+1}-S_{2m+1}^-)^{-1}(S_{2m+2}-S_{2m+2}^-).
\eeao
and that the non symmetric canonical moments defined
by \eqref{nonsym}  satisfy
$$
2\overline{V}_{n-1}\overline{U}_n = (S_{n-1}-S_{n-1}^-)^{-1}(S_n-S_n^-),
$$
whenever $n\leq N(\mu_I)$, where $\overline{V}_n=I_p-\overline{U}_n$.
Consequently, the assertion follows by  a direct application of the second part of
Theorem 4.3.
\hfill $\Box$

\bigskip

 {\bf Acknowledgements.}
The authors are grateful to Martina Stein
who typed parts of this paper with considerable technical expertise.
The work of Holger Dette was supported by the Sonderforschungsbereich TR 12,
(Teilprojekt C2) and in part by a DFG grant DE 502/22-3.

\bigskip

 \bibliographystyle{apalike}

\bibliography{unitcircle}

\end{document}